\documentclass [10pt] {article}

 \usepackage [cp1251] {inputenc}
 \usepackage [english] {babel}
 \usepackage {hyperref}
 \usepackage {amssymb}
 \usepackage {amsmath}
 \usepackage {amsthm}
 \usepackage [curve] {xypic}

 \def \:{\colon}
 \def \le{\leqslant}
 \def \ge{\geqslant}
 \def \<{\langle}
 \def \>{\rangle}
 \def \%{\mathbin{\backslash\negmedspace\backslash}}
 \def \*{\star}
 \def \${\bullet}
 \def \+{\mathrel\#}
 \def \_{\underline}
 \def \0{{\textdegree}}

 \def \Z{{\mathbb Z}}
 \def \k{{\boldsymbol k}}
 \def \U{{\boldsymbol U}}
 \def \V{{\boldsymbol V}}
 \def \E{{\boldsymbol E}}

 \def \kMod{{\text{\rm$\k$-\bf Mod}}}
 \def \kAlc{{\text{\rm$\k$-\bf Alc}}}
 \def \Graph{{\mathbf{Graph}}}

 \let \Im\undefined
 \DeclareMathOperator {\Im} {Im}
 \DeclareMathOperator {\Ob} {Ob}
 \DeclareMathOperator {\Mor} {Mor}
 \DeclareMathOperator {\Fun} {{\bf Fun}}

 \newcommand* {\subhead} [1]
              {\addvspace \bigskipamount
               \noindent
               {\mathversion{bold}\bf\itshape #1\/}}

 \newenvironment* {claim} [1] []
 {\begin{trivlist}\item [\hskip\labelsep {\bf #1}] \it}
 {\end{trivlist} }

 \newenvironment* {demo} [1] []
 {\begin{trivlist}\item [\hskip\labelsep {\it #1}] }
 {\end{trivlist} }

 \usepackage {graphicx}
 \def \figure{\addvspace \bigskipamount
              \noindent
              \includegraphics [viewport = 36 36 559 226,
                                width = \textwidth] {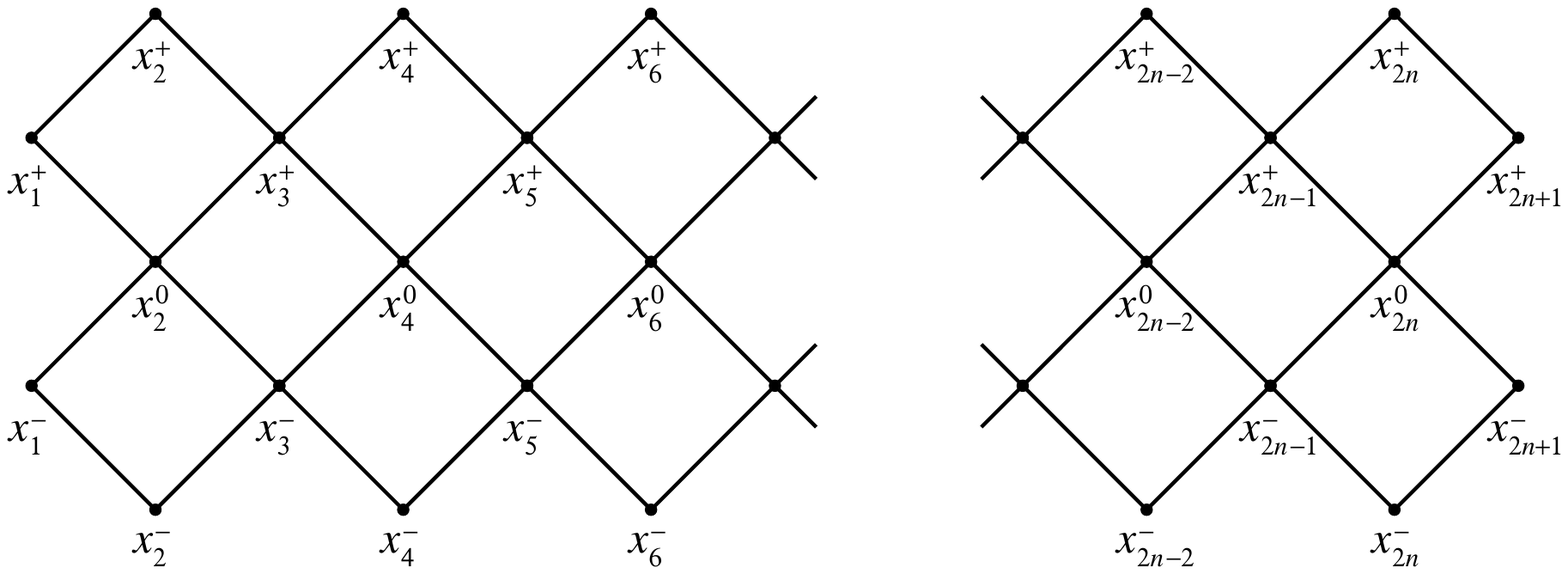}
              \par
              \addvspace \medskipamount}

\begin {document} \frenchspacing

 \title {\Large\bf On the algebra of the M\"obius crown}

 \author {\normalsize\rm S.~S.~Podkorytov}

 \date {}

 \maketitle

 \begin {abstract} \noindent
 A commutative algebra over a field gives rise
 to a representation of the category
 of finite sets and
 surjective maps.
 We consider the restriction of this representation
 to the subcategory of sets of cardinality at most $r$.
 For each $r$,
 we present two non-isomorphic algebras
 that give rise to isomorphic representations of this
 subcategory.
 \end {abstract}


 Let $\Omega_r$
 ($r=0,1,\dotsc,\infty$)
 be the category
 whose objects are the sets $\<p\>=\{1,\dotsc,p\}$,
 $p=1,2,\dotsc$, $p\le r$,
 and
 whose morphisms are surjective maps.
 Let $\k$ be a field.
 We imply it
 when saying about
 vector spaces,
 tensor products,
 etc.
 By an {\it algebra\/} we mean
 a commutative non-unital $\k$-algebra.
 An algebra $A$ gives rise to the functor $L^r(A)\:\Omega_r\to\kMod$
 (a {\it representation\/} of $\Omega_r$)
 that
 takes an object $\<p\>$ to the vector space $A^{\otimes p}$ and
 takes a morphism $s\:\<p\>\to\<q\>$ to the linear map
 $$
 A^{\otimes p}\to A^{\otimes q},
 \qquad
 a_1\otimes\dotso\otimes a_p\mapsto
 m_1\otimes\dotso\otimes m_q,
 $$
 where
 $$
 m_j=\prod_{i\in s^{-1}(j)}a_i
 $$
 (a variant of the Loday functor \cite[Proposition~6.4.4]{Loday}).

 Must algebras $A$ and $B$ be isomorphic if
 the representations $L^r(A)$ and $L^r(B)$ are isomorphic?
 Yes if
 $r=\infty$,
 the field $\k$ is algebraically closed and
 the algebras have finite (vector-space) dimension
 (\cite{Podkorytov}, cf.\ \cite{Dreckmann}).
 Our aim here is to show that
 this is false for arbitrarily large finite $r$.
 For
 each $r=1,2,\dotso$ and
 arbitrary $\k$,
 we present two non-isomorphic finite-dimensional algebras $A$
 and $B$ with isomorphic representations $L^r(A)$ and $L^r(B)$.
 These algebras are obtained from the Stanley--Reisner algebras
 of certain graphs (``crowns'')
 by taking the homogeneous components of degrees 1 and 2.


 \subhead {The functor $L^r$.}
 The correspondence $A\mapsto L^r(A)$ is covariant
 in an obvious way.
 So we have the functor $L^r\:\kAlc\to\Fun(\Omega_r,\kMod)$,
 where
 $\kAlc$ is the category of algebras and
 $\Fun(\Omega_r,\kMod)$ is that of functors $\Omega_r\to\kMod$
 (representations).

 \subhead {The action category $M\%S$.}
 Let a monoid $M$ act on a set $S$ from the left.
 For $s,t\in S$,
 put $M(s,t)=\{m:m\cdot s=t\}\subseteq M$.
 We have the category $M\%S$,
 where
 $\Ob M\%S=S$,
 a bijection
 $$
 M(s,t)\to\Mor_{M\%S}(s,t),
 \qquad
 m\mapsto m|_{s\to t},
 $$
 is given for each $s,t\in S$,
 $1_s=1|_{s\to s}$, and
 the composition of morphisms is given by the multiplication in
 $M$.

 We have the not necessarily commutative unital $\k$-algebra
 $\k[M]$.
 For $s,t\in S$,
 we have the subspace $\k[M(s,t)]\subseteq\k[M]$.

 Consider the linear category $\k[M\%S]$.
 For $s,t\in S$,
 we have the linear map
 $$
 \k[M(s,t)]\to\Mor_{\k[M\%S]}(s,t),
 \qquad
 X\mapsto X\|_{s\to t},
 $$
 given by the rule $[m]\mapsto[m|_{s\to t}]$.
 Clearly,
 $1\|_{s\to s}=1_s$
 ($s\in S$).
 If
 $X\in\k[M(s,t)]$,
 $Y\in\k[M(t,u)]$
 ($s,t,u\in S$),
 then
 $YX\in\k[M(s,u)]$ and
 $$
 (YX)\|_{s\to u}=Y\|_{t\to u}\circ X\|_{s\to t}.
 $$


 \subhead {The monoid $W_n$ and the elements $T_n$ and $Z_n$.}
 Introduce
 the multiplicative submonoid $\V=\{1,-1,0\}\subseteq\Z$ and
 its submonoids
 $\U=\{1,-1\}$ and
 $\E=\{1,0\}$.
 We denote the elements $1$ and $-1$ also by $+$ and $-$
 (respectively).

 Let $W_n\subseteq\V^{2n+1}$ be the submonoid formed by the
 collections
 $$
 w=(w_1,w_2,\dotsc,w_{2n+1})
 $$
 in which
 $w_{2i+1}\in\U$
 ($i=0,\dotsc,n$) and
 $w_jw_{j+1}\in\E$
 ($j=1,\dotsc,2n$).

 Introduce the elements $g_i,h_i\in W_n$ ($i=1,\dotsc,n$):
 \def \I{\vphantom{0^{\frac00}}}
 $$
 g_i=({+},\dotsc,{+},\underset{\I2i}0,{+},\dotsc,{+}),
 \qquad
 h_i=({-},\dotsc,{-},\underset{\I2i}0,{+},\dotsc,{+})
 $$
 and
 $T_n,Z_n\in\k[W_n]$:
 $$
 T_n=\sum_{i=1}^n(1-[g_1])\dotso(1-[g_{i-1}])[h_i],
 \qquad
 Z_n=(1-[g_1])\dotso(1-[g_n]).
 $$
 Using
 commutativity of $W_n$ and
 the relations
 $g_i^2=h_i^2=g_i$ and
 $g_ih_i=h_i$,
 we get
 $$
 T_n^2=1-Z_n.
 $$


 \subhead {Two actions of $W_n$ and their categories.}
 The monoid $W_n$ acts on the set $\U$ from the left by the
 rule $w\cdot s=w_1w_{2n+1}s$.
 Since
 $T_n\in\k[W_n(s,-s)]$ and
 $Z_n\in\k[W_n(s,s)]$
 for each $s\in\U$,
 we have
 \begin {equation} \label {1}
 T_n\|_{-s\to s}\circ T_n\|_{s\to-s}=1_s-Z_n\|_{s\to s}
 \end {equation}
 in $\k[W_n\%\U]$.

 Consider the one-element set $\{\*\}$ with the left action of
 $W_n$.
 The map $\U\to\{\*\}$ induces the functors
 $\omega_n\:W_n\%\U\to W_n\%\{\*\}$ and
 $\k[\omega_n]\:\k[W_n\%\U]\to\k[W_n\%\{\*\}]$.
 For any $s,t\in\U$ and $X\in\k[W_n(s,t)]$,
 we have
 \begin {equation} \label {2}
 \k[\omega_n]\:
 X\|_{s\to t}\mapsto X\|_{\*\to\*}.
 \end {equation}


 \subhead {Graphs.}
 By a {\it graph\/} we mean a pair $G=(G_1,G_2)$,
 where
 $G_1$ is a set and
 $G_2\subseteq G_1\times G_1$ is a reflexive symmetric
 relation.
 The {\it vertices\/} of $G$ are the elements of $G_1$;
 its {\it edges\/} are the sets $\{x,y\}$,
 where $(x,y)\in G_2$, $x\ne y$.

 A morphism $f\:G\to H$ of graphs is a pair $f=(f_1,f_2)$,
 where $f_p\:G_p\to H_p$, $p=1,2$, are maps such that
 $f_2(x,y)=(f_1(x),f_1(y))$, $(x,y)\in G_2$.
 Graphs and their morphisms form a category $\Graph$.


 \subhead {The cofunctor $Q$: the algebra of a graph.}
 Let $G$ be a graph.
 The symmetric group $\Sigma_2$ acts on
 $G_2\subseteq G_1\times G_1$ by permuting the cooordinates.
 We have the projection
 $$
 \k^{G_2}\to(\k^{G_2})_{\Sigma_2},
 \qquad
 u\mapsto\bar u.
 $$
 Let $A^\$$ be the graded algebra
 concentrated in degrees 1 and 2:
 $$
 A^1=\k^{G_1},
 \qquad
 A^2=(\k^{G_2})_{\Sigma_2},
 $$
 where,
 if $a,b\in A^1$,
 then $ab=\bar u\in A^2$,
 where $u\in\k^{G_2}$, $u(x,y)=a(x)b(y)$.

 Put $Q^\$(G)=A^\$$.
 Let $Q(G)$ be the same algebra considered without the grading.
 The correspondence $G\mapsto Q(G)$ is contravariant
 in an obvious way.
 So we have the cofunctor $Q\:\Graph\to\kAlc$.
 We need the following properties of $Q$.

 \hypertarget{P1}{1\0}.
 If $G$ is finite,
 then $Q(G)$ has finite dimension.

 \hypertarget{P2}{2\0}.
 If graph morphisms $f_i\:G_i\to H$, $i\in I$, form a {\it
 cover},
 i.~e.,
 $$
 \bigcup_{i\in I}\Im f_{i\;p}=H_p,
 \qquad
 p=1,2,
 $$
 then the linear map
 $$
 (Q(f_i))_{i\in I}\:Q(H)\to\prod_{i\in I}Q(G_i)
 $$
 is injective.

 \hypertarget{P3}{3\0}.
 If finite graphs $G$ and $H$ are non-isomorphic,
 then the algebras $Q(G)$ and $Q(H)$ are non-isomorphic too.
 This follows from the Gubeladze theorem
 \cite[Theorem~3.1]{Gubeladze}.
 We give simpler arguments
 that suffice in the special case
 that we will need.

 Call a graph $G$ {\it admissible\/} if,
 for any distinct $x,y\in G_1$,
 there exists $z\in G_1$ such that
 $(x,z)\notin G_2$ and
 $(y,z)\in G_2$.
 (For example,
 any graph without triangles and pendant vertices is
 admissible.)
 We show that
 an admissible graph $G$ can be reconstructed from $Q(G)$.

 Let $A^\$$ be a graded algebra
 concentrated in degrees 1 and 2.
 Consider the projective space $P(A^1)$.
 Let $[\ ]\:A^1\setminus\{0\}\to P(A^1)$ be the projection.
 Define on $P(A^1)$
 a symmetric relation $\+$ ({\it dependence\/}):
 $[a]\+[b]\Leftrightarrow ab\ne0$,
 and
 a preorder $\lesssim$:
 $p\lesssim q\Leftrightarrow p^{\+}\subseteq q^{\+}$,
 where $r^{\+}=\{s:r\+s\}$.
 Let $R\subseteq P(A^1)$ be the set of {\it minimal\/} points,
 i.~e.\ those points $p$ for which
 $\{s:s\lesssim p\}=\{p\}$.
 If $A^\$=Q^\$(G)$ for some graph $G$,
 then there is the injective map
 $e\:G_1\to P(A^1)$, $x\mapsto[\delta_x]$,
 where
 $\delta_x\in A^1=\k^{G_1}$, $\delta_x(y)$ equals
 $1$ if $y=x$ and
 $0$ otherwise.
 The inverse image of $\+$ under $e$ equals $G_2$.
 It is not hard to check that,
 if $G$ is admissible,
 then $\Im e=R$.
 It remains to add that
 the graded algebra $A^\$$ can be reconstructed from the
 ungraded algebra $A=Q(G)$:
 $A^\$$ is canonically isomorphic to the graded algebra $B^\$$
 with
 the components $B^1$ and $B^2$,
 where
 $B^2=\{b:bA=0\}\subseteq A$ and
 $B^1=A/B^2$
 (so
 $B^2=A^2$ and
 $B^1\cong A^1$),
 and
 the multiplication induced by that in $A$.


 \subhead {The graph $B_n$.}
 Let $B_n$ be the graph shown on the figure.
 Its vertices are $x_j^v$,
 where
 $j=1,\dotsc,2n+1$,
 $v\in\V$,
 and
 $v\in\U$
 if $j$ is odd.

 \figure

 The monoid $W_n$ acts on $B_n$ from the left by the rule
 $w\cdot x_j^v=x_j^{w_jv}$.
 Let $w_*\:B_n\to B_n$ be the action of $w\in W_n$.
 The graph $B_n$ with the action of $W_n$ gives rise to the
 functor
 $$
 \_B_n\:W_n\%\{\*\}\to\Graph,
 \qquad
 \*\mapsto B_n,
 \qquad
 w|_{\*\to\*}\mapsto w_*.
 $$
 Since $\Fun(\Omega_r,\kMod)$ is a linear category,
 the cofunctor
 $$
 W_n\%\{\*\}
 \xrightarrow{\_B_n}
 \Graph
 \xrightarrow Q
 \kAlc
 \xrightarrow{L^r}
 \Fun(\Omega_r,\kMod)
 $$
 extends to a linear cofunctor
 $$
 b_n^r\:\k[W_n\%\{\*\}]\to\Fun(\Omega_r,\kMod).
 $$

 \begin {claim} [\hypertarget{Lemma}{Lemma}.]
 We have $b_n^{n-1}(Z_n\|_{\*\to\*})=0$.
 \end {claim}
 
 \begin {demo} [Proof.]
 Take $p=1,\dotsc,n-1$.
 The monoid $W_n$ acts on $B_n$ from the left.
 The induced right action on the vector space
 $Q(B_n)^{\otimes p}$ makes it a right $\k[W_n]$-module.
 We should show that
 $Q(B_n)^{\otimes p}Z_n=0$.

 For $i=1,\dotsc,n$,
 let $F_i$ be the subgraph of $B_n$
 spanned by the vertices $x_j^v$ with $|j-2i|\le1$ and
 let $e_i\:F_i\to B_n$ be the inclusion morphism.
 Since the subgraphs $F_i$ cover $B_n$,
 the linear map
 $$
 (Q(e_i))_{i=1}^n\:Q(B_n)\to\bigoplus_{i=1}^nQ(F_i)
 $$
 is injective
 (by the property~\hyperlink{P2}{2\0}).
 Raising it to the tensor power $p$,
 we get an injective linear map
 $$
 E_p\:Q(B_n)^{\otimes p}\to
 \bigoplus_{i_1,\dotsc,i_p}S_{i_1\dotso i_p},
 \qquad
 S_{i_1\dotso i_p}=Q(F_{i_1})\otimes\dotso\otimes Q(F_{i_p}).
 $$
 The subgraphs $F_i$ are invariant under the action of $W_n$.
 The induced right action on the vector spaces
 $S_{i_1\dotso i_p}$ makes them right $\k[W_n]$-modules.
 The map $E_p$ is a homomorphism of $\k[W_n]$-modules.
 Since it is injective,
 it suffices to show that
 $S_{i_1\dotso i_p}Z_n=0$.

 Each element $g_i$ acts trivially on the subgraphs $F_{i'}$,
 $i'\ne i$.
 Thus,
 if $i$ is distinct from $i_1,\dotsc,i_p$,
 the element $g_i$ acts trivially on $S_{i_1\dotso i_p}$ and
 thus $S_{i_1\dotso i_p}Z_n=0$.
 Since $p<n$,
 such an $i$ exists for any $i_1,\dotsc,i_p$.
 \qed
 \end {demo}


 \subhead {The graphs $C_n^s$ (crowns).}
 Take $n\ge2$.
 For $s\in\U$,
 let $C_n^s$ be the graph
 obtained from $B_n$ by identifying
 $x_{2n+1}^v$ with $x_1^{sv}$ for each $v\in\U$.
 Let $f_n^s\:B_n\to C_n^s$ be the projection morphism.
 We call
 $C_n^+$ the {\it simple crown\/} and
 $C_n^-$ the {\it M\"obius\/} one.

 The graphs $C_n^s$, $s\in\U$, are non-isomorphic
 (the edges containing vertices of valency 2
 form
 two cycles in $C_n^+$ and
 one cycle in $C_n^-$).
 They are finite and admissible, and
 thus
 (see the properties
 \hyperlink{P1}{1\0} and
 \hyperlink{P3}{3\0})
 their algebras $Q(C_n^s)$ are finite-dimensional and
 non-isomorphic.
 We show that
 the representations $L^{n-1}(Q(C_n^s))$, $s\in\U$, are
 isomorphic.

 For $s,t\in\U$ and $w\in W_n(s,t)$,
 let $w_*\:C_n^s\to C_n^t$ be the morphism such that
 the following diagram is commutative:
 $$
 \xymatrix {
 B_n
 \ar[d]_-{w_*}
 \ar[rr]^-{f_n^s} &&
 C_n^s
 \ar[d]^-{w_*} \\
 B_n
 \ar[rr]^-{f_n^t} &&
 C_n^t.
 }
 $$
 So we have the functor
 $$
 \_C_n\:W_n\%\U\to\Graph,
 \qquad
 s\mapsto C_n^s,
 \qquad
 w|_{s\to t}\mapsto w_*.
 $$
 The morphisms $f_n^s$, $s\in\U$, form a morphism of functors
 $f_n\:\_B_n\circ\omega_n\to\_C_n$:
 $$
 \xymatrix {
 &&
 \Graph
 && \\
 &&&& \\
 W_n\%\{\*\}
 \ar[uurr]^-{\_B_n}
 \ar@=[urrr]^-{f_n} &&&&
 W_n\%\U.
 \ar[llll]^-{\omega_n}
 \ar[uull]_-{\_C_n}
 }
 $$

 Since $\Fun(\Omega_r,\kMod)$ is a linear category,
 the cofunctor
 $$
 W_n\%\U
 \xrightarrow{\_C_n}
 \Graph
 \xrightarrow Q
 \kAlc
 \xrightarrow{L^r}
 \Fun(\Omega_r,\kMod)
 $$
 extends to a linear cofunctor
 $$
 c_n^r\:\k[W_n\%\U]\to\Fun(\Omega_r,\kMod).
 $$
 The morphism $f_n$ induces a morphism of cofunctors
 $$
 \xymatrix {
 &&
 \Fun(\Omega_r,\kMod)
 && \\
 &&&
 {\phantom{OO}}
 \ar@=[dlll]
 & \\
 \k[W_n\%\{\*\}]
 \ar[uurr]^-{b_n^r} &&&&
 \k[W_n\%\U],
 \ar[llll]^-{\k[\omega_n]}
 \ar[uull]_-{c_n^r}
 }
 $$
 i.~e.,
 for any $s,t\in\U$ and $X\in\k[W_n(s,t)]$,
 we have the commutative diagram
 $$
 \xymatrix {
 L^r(Q(B_n)) &&
 L^r(Q(C_n^s))
 \ar[ll]_-{L^r(Q(f_n^s))} \\
 L^r(Q(B_n))
 \ar[u]^-{b_n^r(X\|_{\*\to\*})} &&
 L^r(Q(C_n^t))
 \ar[u]_-{c_n^r(X\|_{s\to t})}
 \ar[ll]_-{L^r(Q(f_n^t))}
 }
 $$
 (we used the rule~\eqref{2}).
 Since $f_n^s$ is a cover,
 the homomorphism $Q(f_n^s)\:Q(C_n^s)\to Q(B_n)$ is injective
 (by the property~\hyperlink{P2}{2\0}),
 and
 thus
 the morphism $L^r(Q(f_n^s))$ is objectwise injective.

 Now assume $r=n-1$, $s=t$ and $X=Z_n$.
 By \hyperlink{Lemma}{Lemma},
 $b_n^{n-1}(Z_n\|_{\*\to\*})=0$.
 Thus
 $c_n^{n-1}(Z_n\|_{s\to s})=0$
 (by
 commutativity of the diagram and
 the mentioned objectwise injectivity).
 We show that
 the arrows of the diagram
 $$
 \xymatrix {
 L^{n-1}(Q(C_n^+))
 \ar@/^2ex/[rr]^-{c_n^{n-1}(T_n\|_{-\to+})} &&
 L^{n-1}(Q(C_n^-))
 \ar@/^2ex/[ll]^-{c_n^{n-1}(T_n\|_{+\to-})}
 }
 $$
 are mutually inverse.
 For each $s\in\U$,
 we have
 \begin {multline*}
 c_n^{n-1}(T_n\|_{s\to-s})\circ c_n^{n-1}(T_n\|_{-s\to s})=
 c_n^{n-1}(T_n\|_{-s\to s}\circ T_n\|_{s\to-s})= \\ =
 c_n^{n-1}(1_s-Z_n\|_{s\to s})=
 1_{L^{n-1}(Q(C_n^s))}-c_n^{n-1}(Z_n\|_{s\to s})=
 1_{L^{n-1}(Q(C_n^s))}
 \end {multline*}
 (we used the equality~\eqref{1}).


 \bigskip

 \noindent
 I am grateful to I.~S.~Baskov,
 my conversations with whom resulted in this work.


 \bigskip

 \noindent
 St.~Petersburg Department of Steklov Mathematical Institute
 \newline
 of Russian Academy of Sciences

 \smallskip

 \noindent \href {mailto:ssp@pdmi.ras.ru} {\tt ssp@pdmi.ras.ru}

 \noindent \url {http://www.pdmi.ras.ru/~ssp}

 \end {document}